\documentclass[number,citesort,seceqn,dvips]{arxbj}
\usepackage{dcolumn}
\usepackage[left]{vmkcol}

\aid{0}
\volume{18}
\issue{1}
\pubyear{2012}
\firstpage{206}
\lastpage{228}
\doi{10.3150/10-BEJ330}

\makeatletter
\newcolumntype{d}[1]{D{.}{.}{#1}}

\newtheorem{theorem}{Theorem}[section]
\newtheorem{lem}{Lemma}[section]
\newremark{rem}{Remark}[section]
\makeatother

\begin{document}
\begin{frontmatter}

\title{Efficient estimation of moments in linear mixed models}
\runtitle{Linear mixed models}

\begin{aug}
\author[1]{\fnms{Ping} \snm{Wu}\ead[label=e1]{wu\_ping916@yahoo.com.cn}\thanksref{1}},
\author[2]{\fnms{Winfried} \snm{Stute}\corref{}\thanksref{2}\ead[label=e2,text={winfried.stute@math.uni-giessen.de}]{winfried.stute@math.uni-giessen.de}}\and
\author[3]{\fnms{Li-Xing} \snm{Zhu}\ead[label=e3]{lzhu@hkbu.edu.hk}\thanksref{3}}

\runauthor{P. Wu, W. Stute and L.-X. Zhu}

\address[1]{East China Normal University, Shanghai, China. \printead{e1}}
\address[2]{Mathematical Institute, University of Giessen, Arndtstr. 2,
D-35392 Giessen, Germany.\\ \printead{e2}}
\address[3]{Hong Kong Baptist University, Hong Kong, China. \printead{e3}}
\end{aug}

% HISTORY:
\received{\smonth{11} \syear{2009}}
\revised{\smonth{8} \syear{2010}}

% ABSTRACT
%
\begin{abstract}
In the linear random effects model, when distributional assumptions
such as normality of the error variables cannot be justified, moments
may serve as alternatives to describe relevant distributions in
neighborhoods of their means. Generally, estimators may be obtained as
solutions of estimating equations. It turns out that there may be
several equations, each of them leading to consistent estimators, in
which case finding the efficient estimator becomes a crucial problem.
In this paper, we systematically study estimation of moments of the
errors and random effects in linear mixed models.
\end{abstract}

% KEYWORDS
%

\end{frontmatter}

%s1 ###
\section{Introduction}

Normality or, more generally, the existence of a parametric structure on
the distribution of random effects is a routine assumption for linear
mixed models. In such a case, both the maximum likelihood estimator
(MLE) and the restricted maximum likelihood estimator (RMLE) work well.
Moreover, they are standard outputs in statistical software packages
such as SAS and R. A comprehensive account of the methodology is
contained in the monograph of Verbeke and Molenberghs \cite{vermol}.
In recent years, more efforts were devoted to relaxing this assumption
and using semiparametric or nonparametric methods to estimate the
parameters of interest. Zhang and Davidian \cite{zhadav} suggested
using the seminonparametric representation of Gallant and Nychka
\cite{gallnyc} to approximate the random effect density in order to
estimate parameters for linear mixed models. Cui, Ng and Zhu~\cite{cuingzhu} used the estimation of moments in mixed effect models
with errors in variables. Rank estimation was applied by Wang and Zhu
\cite{wanzhu}
to estimate fixed effects.

However, the aforementioned papers do not consider the estimation of
higher moments that are useful for hypothesis testing and interval
estimation for the parameters in the models. To the best of our
knowledge, Cox and Hall \cite{coxhall} is the only reference in the
literature that\vadjust{\goodbreak}
defines and studies the estimators of the errors and random effects
for higher than second moments. The authors of that work obtained the
cumulants of the two components of variance based on homogeneous
polynomials in a simple random effects location model, which is the
sum of the one-level random effect and the error. For this model,
Hall and Yao \cite{hallyao} studied nonparametric estimation of the
distributions
of the errors and the random effects via empirical cumulant generating
functions. To the best of our knowledge, no paper has investigated
this issue for the linear mixed model under consideration.

The contents of this paper are as follows:
\begin{itemize}
\item
In Section \ref{sec21} we introduce the linear mixed model and derive basic
properties of the generalized least squares estimator under weak
conditions on the group sizes and the design variables. The fundamental
Lemma \ref{lem21} yields representations of certain polynomial
functions of the
overall errors in terms of individual and group errors. This will be
the basic tool to answer a question posed by Cox and Hall \cite
{coxhall} in the
context of the simple random effects location model, namely, how to
properly weight and combine certain polynomial functions of the residuals.
\item
As a warmup, in Section \ref{sec22}, we consider the estimation of second
moments. It turns out that by a proper combination of polynomial
functions of the residuals, we can obtain second moment estimators
which are asymptotically normal and have the same limit variance as if
the unknown errors were known.
\item
For third and fourth moments, the situation is more complex. In
Sections \ref{sec23} and~\ref{sec24}, we propose and study estimators
yielding efficiency and asymptotic normality under weak conditions
on the design and group sizes.
\item
As an alternative, in Sections \ref{sec31} and \ref{sec32}, we study
an extension of
an estimator due to Cox and Hall \cite{coxhall} which may therefore be
considered as a first step estimator. When the group sizes are all
equal, our estimators have similar asymptotic properties to theirs. We
show that for unequal group sizes, the obtained estimators may converge
at slower rates unless some restrictive regularity assumptions are satisfied.
\item
Section \ref{sec4} presents some simulation studies, while proofs are deferred
to the \hyperref[appendix]{Appendix}.
\end{itemize}

%s2 ###
\section{Minimum variance estimation of moments}\label{sec2}

%s2.1 ###
\subsection{Motivation and first results}\label{sec21}

Assume that data are available from a linear mixed model, that is, we
observe pairs $(x_{ij}, y_{ij}), 1\le i\le n, 1 \le j \le l_i$,
satisfying
%
%e1 ###
%
\begin{equation}\label{eq2.1}
y_{ij} = \alpha+ x'_{ij} \beta+ b_i + \varepsilon_{ij} .
\end{equation}
Here, $i$ denotes the group index, while the measurements within this
group are indexed by~$j$. The integer $l_i$ is the sample size within
group $i$. The row vector $x'_{ij}$ is a $p$-dimensional input vector
corresponding to the $j$th observation in the $i$th group leading to
the output~$y_{ij}$. The relation between $x_{ij}$ and $y_{ij}$
described by (\ref{eq2.1}) contains the intercept parameter $\alpha$,
the fixed effect regression parameter $\beta$\vadjust{\goodbreak} and the one-level random
effect $b_i$ for group $i$, all unknown. Moreover, these quantities are
disturbed by random errors $\varepsilon_{ij}$. It is assumed throughout
that $b_1,\ldots, b_n$ are independent and identically distributed
(i.i.d.) and also independent of all $\varepsilon_{ij}$, which are also
i.i.d. Finally, we may assume without loss of generality that\vspace*{-1pt}
%
%e2 ###
%
\begin{equation}\label{eq2.2}
\mathbb{E} b_i = 0\quad\mbox{and}\quad\mathbb{E} \varepsilon_{ij}
= 0\qquad\mbox{for } 1\le i \le n \mbox{ and } 1\le j \le l_i .\vspace*{-1pt}
\end{equation}
Otherwise,
we may incorporate unknown nonzero expectations in the intercept
$\alpha$. Let~$\gamma_b^k$ and $\gamma_{\varepsilon}^k$ denote the
$k$th moments of the random effects and errors, respectively. In this
paper, we shall construct and analyze estimators of $\alpha$, $\beta$,
$\gamma_b^k$ and $\gamma_{\varepsilon}^k$, $k = 2,3,4,$ that are based
on various estimating equations. These equations are obtained from
proper nonlinear combinations of the residuals. For these, we first
have to estimate $\beta$ and $\alpha$ via a generalized least squares
method. In the model (\ref{eq2.1}), this leads to\vspace*{-1pt}
%
%e3 ###
%
\begin{equation}\label{eq2.3}
\hat\beta= \hat\Sigma_n^{-1} \frac{\sum_{i=1}^n \sum
_{j=1}^{l_i} (x_{ij} - \bar{x}_{i\cdot}) (y_{ij} - \bar{y}_{i\cdot
})}{\sum_{i=1}^n l_i}\vspace*{-1pt}
\end{equation}
and\vspace*{-1pt}
%
%e4 ###
%
\begin{equation}\label{eq2.4}
\hat\alpha= \frac{1}{n} \sum_{i=1}^n \bar{y}_{i\cdot} - \frac{1}{n}
\sum_{i=1}^n \bar{x}'_{i\cdot} \hat\beta.\vspace*{-1pt}
\end{equation}
Here,\vspace*{-1pt}
%
%e5 ###
%
\begin{equation}\label{eq2.5}
\hat\Sigma_n = \frac{1}{\sum_{i=1}^n l_i} \sum_{i=1}^n \sum
_{j=1}^{l_i} (x_{ij} - \bar{x}_{i\cdot})(x_{ij} - \bar{x}_{i\cdot})',\vspace*{-1pt}
\end{equation}
while\vspace*{-1pt}
\[
\bar{x}_{i\cdot} = \frac{1}{l_i} \sum_{j=1}^{l_i} x_{ij}\quad\mbox
{and}\quad
\bar{y}_{i\cdot} = \frac{1}{l_i} \sum_{j=1}^{l_i} y_{ij}\vspace*{-1pt}
\]
denote the corresponding group averages. Furthermore, we let\vspace*{-1pt}
\[
N = \sum_{i=1}^n l_i,\vspace*{-1pt}
\]
the overall sample size.

\begin{theorem}\label{thm:1}
Assume that the following conditions
(\ref{eq2.6})--(\ref{eq2.8}) are satisfied:\vspace*{-1pt}
%
%e6 ###
%
\begin{equation}\label{eq2.6}
\lim_{n\to\infty}\hat{\Sigma}_n = \Sigma \mbox{ for
some positive definite $p\times p$ matrix $\Sigma$};\vspace*{-6pt}
\end{equation}
%
%e7 ###
%
\begin{equation}\label{eq2.7}
\frac{1}{n}\sum_{i=1}^n \frac{1}{l_i} \to0\quad\mbox{and}\quad
N/n \to
\infty\qquad\mbox{as } n \to\infty;\vspace*{-6pt}
\end{equation}
%
%e8 ###
%
\begin{equation}\label{eq2.8}
\frac{\max_{i\le i\le n, 1\le j\le l_i} \| x_{ij} - \bar{x}_{i\cdot
}\|
}{\sqrt{N}} \to0\quad\mbox{and}\quad\frac{1}{n} \sum_{i=1}^n
\bar{x}_{i\cdot
} \mbox{ is bounded}.\vspace*{-1pt}\vadjust{\goodbreak}
\end{equation}
\end{theorem}

Then, in distribution, we have\vspace*{-1pt}
%
%e9 ###
%
\begin{equation}\label{eq2.9}
N^{1/2} (\hat\beta- \beta) \to{\mathcal N}_p(0, \gamma_{\varepsilon}^2
\Sigma^{-1})\vspace*{-1pt}
\end{equation}
and\vspace*{-1pt}
%
%e10 ###
%
\begin{equation}\label{eq2.10}
n^{1/2} (\hat\alpha- \alpha) \to{\mathcal N}_1(0,\gamma_b^2) .\vspace*{-1pt}
\end{equation}
The estimators $\hat\beta$ and $\hat\alpha$ and their distributional
behavior play an important role for motivating the estimation of
$\gamma_b^k$ and $\gamma_{\varepsilon}^k$ since this will be based on
the residuals\vspace*{-1pt}
\[
\hat e_{ij} = y_{ij} - \hat\alpha- x'_{ij}\hat\beta.\vspace*{-1pt}
\]
Set\vspace*{-1pt}
\[
\bar{b} \equiv\bar{b}_n= \frac{1}{n} \sum_{i=1}^n b_i,\qquad \bar
{\varepsilon} \equiv\bar{\varepsilon}_n = \frac{1}{n} \sum_{i=1}^n
\frac{1}{l_i} \sum_{j=1}^{l_i} \varepsilon_{ij}\vspace*{-1pt}
\]
and\vspace*{-1pt}
\[
\bar{x} \equiv\bar{x}_n = \frac{1}{n} \sum_{i=1}^n \frac{1}{l_i}
\sum
_{j=1}^{l_i} x_{ij} .\vspace*{-1pt}
\]
In view of (\ref{eq2.4}), we have\vspace*{-1pt}
\[
\hat\alpha- \alpha= \bar{b} + \bar{\varepsilon} - \bar{x}'(\hat
\beta
- \beta),\vspace*{-1pt}
\]
from which it follows that\vspace*{-1pt}
\begin{eqnarray} \label{eq2.11}
\hat e_{ij} &=& (b_i - \bar{b}) + (\varepsilon_{ij} - \bar
{\varepsilon
}) + (x_{ij} - \bar{x})'(\beta- \hat\beta) \nonumber\\ [-9pt]\\ [-9pt]
&\equiv&(b_i + \varepsilon_{ij}) - (\bar{b} + \bar{\varepsilon}) +
z_{ij}'(\beta- \hat\beta).\nonumber\vspace*{-1pt}
\end{eqnarray}
Set\vspace*{-1pt}
\[
e_{ij} = b_i + \varepsilon_{ij} ,\vspace*{-1pt}
\]
a sum of two independent zero-mean random variables.

When the $l_i$'s are equal and $\beta= 0$, that is, in the simple
random effects location model, Cox and Hall \cite{coxhall} used homogeneous
polynomial functions to construct estimating equations. In the present
paper, we consider more general situations in which new special
nonlinear functions of the $e_{ij}$'s are important tools to derive
estimating equations for\vspace*{-1pt}
\[
\gamma_b^k = \mathbb{E} b_i^k\quad \mbox{and}\quad \gamma_{\varepsilon}^k =
\mathbb{E} \varepsilon_{ij}^k .\vspace*{-1pt}
\]
For this, define, for $1\le i\le n$ and $1 \le m \le k$,\vspace*{-1pt}
\[
f_m^k(i) = \sum_{j=1}^{l_i} e_{ij}^m \Biggl[ \sum_{j=1}^{l_i} e_{ij}\Biggr]^{k-m} .\vspace*{-1pt}
\]
The following lemma turns out to be crucial for our analysis.\vadjust{\goodbreak}

\begin{lem}\label{lem21}
We have\vspace*{-1pt}
\[
f_m^k(i) = \sum_{t=0}^k \sum_{s=(t-k+m)\vee0}^{t\wedge m} {m \choose s}
{k-m \choose t-s} \Biggl( \sum_{j=1}^{l_i} \varepsilon_{ij}^s\Biggr) \Biggl(
\sum_{j=1}^{l_i} \varepsilon_{ij}\Biggr)^{t-s} b_i^{k-t} l_i^{k-m-t+s} .\vspace*{-1pt}
\]
Here, $a \wedge b$ and $a \vee b$ denote the minimum and maximum,
repectively, of two real numbers~$a$ and~$b$.\vspace*{-1pt}
\end{lem}

The proof follows from simple arithmetic. When we take expectations,
usually many of the terms in the expansion of $f_m^k(i)$ will vanish,
mainly because the $\varepsilon_{ij}$'s and $b_i$'s are centered and
independent; see (\ref{eq2.2}). Moreover, by taking proper linear combinations
of the $f_m^k(i)$'s, we shall be able to represent the $\gamma_b^k$'s
and $\gamma_{\varepsilon}^k$'s in terms of the $f$'s. These so-called
estimating equations will then lead to associated estimators.

For example, in the case of $\gamma_{\varepsilon}^2$, we have\vspace*{-1pt}
\[
l_i f_2^2(i) - f_1^2(i) = l_i \sum_{j=1}^{l_i} \varepsilon_{ij}^2 - \Biggl[
\sum_{j=1}^{l_i} \varepsilon_{ij}\Biggr]^2,\vspace*{-1pt}
\]
from which it follows that\vspace*{-1pt}
\[
\mathbb{E} [ l_i f_2^2(i) - f_1^2(i)] = l_i(l_i - 1) \gamma
_{\varepsilon}^2.\vspace*{-1pt}
\]
This equation does not incorporate any $b$-term, so it may serve as a
basis for the estimation of $\gamma_{\varepsilon}^2$. For moments
$\gamma_{\varepsilon}^k$ and $\gamma_b^k$, $k > 2$, things become more
delicate. At first, it is not clear how to combine the $f_m^k(i)$'s in
order to get efficient estimators. This issue is dealt with
in Sections~\ref{sec22}--\ref{sec24}, for $k = 2,3$ and 4, respectively. In Section \ref{sec3}, we briefly
discuss the extension of Cox and Hall \cite{coxhall} to the regression case and
show that it may cause some inefficiencies.\vspace*{-1pt}

\begin{rem}
We only remark in passing that the results of this and the following
sections may be extended to group sizes $l_{ni}$, $1\le i\le n$, that
is, when the $l$'s depend on the number $n$ of groups and therefore
form a triangular array.\vspace*{-3pt}
\end{rem}

%s2.2 ###
\subsection{\texorpdfstring{Estimation of $\gamma^2_{\varepsilon}$ and $\gamma_{b}^2$}
{Estimation of gamma 2 epsilon and gamma b 2}}
\vspace*{-3pt}\label{sec22}

We start by estimating $\gamma_{\varepsilon}^2$ and $\gamma_b^2$. As
mentioned above,\vspace*{-1pt}
\[
\mathbb{E} [ l_i f_2^2(i) - f_1^2(i)] = l_i(l_i - 1) \gamma
_{\varepsilon}^2 .\vspace*{-1pt}
\]
Averaging over $1\le i\le n$ and replacing the unknown $\varepsilon$'s
by the residuals leads to the estimator\vspace*{-1pt}
\[
\hat\gamma_{\varepsilon}^2 = \frac{\sum_{i=1}^n (1/(l_i - 1))
\{
l_i \sum_{j=1}^{l_i} \hat e_{ij}^2 - ( \sum_{j=1}^{l_i} \hat
e_{ij})^2\}}{N} .\vspace*{-1pt}
\]
Similarly, the equation\vspace*{-1pt}
\[
\mathbb{E} [f_1^2(i) - f_2^2(i)] = l_i (l_i - 1) \gamma_b^2\vspace*{-1pt}\vadjust{\goodbreak}
\]
leads to the estimator
\[
\hat\gamma_b^2 = \frac{1}{n} \sum_{i=1}^n \frac{1}{l_i(l_i-1)} \Biggl\{ \Biggl(
\sum_{j=1}^{l_i} \hat e_{ij}\Biggr)^2 - \sum_{j=1}^{l_i} \hat e_{ij}^2 \Biggr\} .
\]

\begin{theorem}\label{thm:2}
Under the conditions of Theorem \ref{thm:1}, when
$\gamma_{\varepsilon}^4$ and $\gamma_b^4$ are finite, we have
that
%
%e11 ###
%
\begin{equation} \label{eq2.12}
N^{1/2} [ \hat\gamma_{\varepsilon}^2 - \gamma_{\varepsilon}^2]
\to{\mathcal N}_1(0,\mu_{\varepsilon}^2)
\end{equation}
and
%
%e12 ###
%
\begin{equation}\label{eq2.13}
n^{1/2} [ \hat\gamma_b^2 - \gamma_b^2] \to{\mathcal N}_1(0,
\mu_b^2) ,
\end{equation}
where
\[
\mu_{\varepsilon}^2 = \gamma_{\varepsilon}^4 - (\gamma
_{\varepsilon
}^2)^2\quad \mbox{and}\quad \mu_b^2 = \gamma_b^4 - (\gamma_b^2)^2 .
\]
\end{theorem}

It is interesting to note that (\ref{eq2.12}) and (\ref{eq2.13}) will
be shown by verifying
\[
N^{1/2} [ \hat\gamma_{\varepsilon}^2 - \gamma_{\varepsilon}^2 ] =
N^{-1/2} \Biggl[ \sum_{i=1}^{n} \sum_{j=1}^{l_i} (\varepsilon_{ij}^2 -
\gamma
_{\varepsilon}^2)\Biggr] + \mathrm{o}_{\mathbb{P}}(1)
\]
and
\[
n^{1/2} [ \hat\gamma_b^2 - \gamma_b^2] = \frac{1}{\sqrt{n}} \Biggl[
\sum_{i=1}^n (b_i^2 - \gamma_b^2)\Biggr] + \mathrm{o}_{\mathbb{P}}(1) .
\]
In other words, $\hat\gamma_{\varepsilon}^2$ and $\hat\gamma_b^2$ are as efficient as the moment estimators based on the true (but
unknown) $\varepsilon_{ij}$ and $b_i$.

%s2.3 ###
\vspace*{2pt}
\subsection{\texorpdfstring{Estimation of $\gamma_{\varepsilon}^3$ and $\gamma_b^3$}
{Estimation of gamma epsilon 3 and gamma b 3}}
\vspace*{2pt}
\label{sec23}

In this section we show how to estimate $\gamma_{\varepsilon}^3$ and
$\gamma_b^3$ with minimal variance. Again, this may be achieved by
properly combining the $f_m^k(i)$'s. From Lemma \ref{lem21}, we obtain
\begin{eqnarray*}
 \mathbb{E} f_3^3(i) &=& l_i \gamma_b^3 + l_i \gamma_{\varepsilon
}^3,\\[2pt]
\mathbb{E} f_2^3(i) &=& l_i^2 \gamma_b^3 + l_i \gamma_{\varepsilon
}^3
\end{eqnarray*}
and
\[
\mathbb{E} f_1^3(i) = l_i^3 \gamma_b^3 + l_i \gamma_{\varepsilon
}^3 .
\]
We conclude that
\[
\mathbb{E} [ 2 f_1^3(i) + l_i^2 f_3^3(i) - 3 l_i f_2^3(i)] = l_i(l_i -
1) (l_i - 2) \gamma_{\varepsilon}^3 .
\]
The corresponding estimator of $\gamma_{\varepsilon}^3$ becomes
\begin{eqnarray*}
\hat\gamma_{\varepsilon}^3 = N^{-1} \sum_{i=1}^n \frac
{1}{(l_i-1)(l_i-2)} \Biggl\{ 2 \Biggl( \sum_{j=1}^{l_i} \hat
e_{ij}\Biggr)^3 + l_i^2 \sum_{j=1}^{l_i} \hat e_{ij}^3 -3 l_i\Biggl( \sum_{j=1}^{l_i} \hat e_{ij}^2\Biggr) \Biggl( \sum_{j=1}^{l_i} \hat
e_{ij}\Biggr)\Biggr\} .
\end{eqnarray*}
For $\gamma_b^3$, the relevant equation is
\[
\mathbb{E} [ f_1^3(i) - 3f_2^3(i) + 2f_3^3(i)] =
l_i(l_i-1)(l_i-2)\gamma_b^3,
\]
leading to the estimator
\[
\hat\gamma_b^3 = \frac{1}{n} \sum_{i=1}^n \frac{1}{l_i(l_i-1)(l_i-2)}
\Biggl\{ \Biggl( \sum_{j=1}^{l_i} \hat e_{ij}\Biggr)^3 - 3\Biggl( \sum_{j=1}^{l_i} \hat
e_{ij}^2\Biggr) \Biggl( \sum_{j=1}^{l_i} \hat e_{ij}\Biggr) + 2 \sum_{j=1}^{l_i} \hat
e_{ij}^3\Biggr\}.
\]

\begin{theorem}\label{thm:5}
 Under the conditions of Theorem \ref{thm:1}, when
$\gamma_{\varepsilon}^6$ and $\gamma_b^6$ are finite, we have that
\[
N^{1/2} (\hat\gamma_{\varepsilon}^3 - \gamma_{\varepsilon}^3) \to
{\mathcal
N}_1(0, \mu_{\varepsilon}^3)
\]
and
\[
n^{1/2} (\hat\gamma_b^3 - \gamma_b^3) \to{\mathcal N}_1(0,\mu_b^3) ,
\]
where
\begin{eqnarray*}
\mu_{\varepsilon}^3 &=& \gamma_{\varepsilon}^6 - (\gamma
_{\varepsilon
}^3)^2 - 6 \gamma_{\varepsilon}^2 \gamma_{\varepsilon}^4 + 9(\gamma
_{\varepsilon}^2)^3,\\
\mu_{b}^3 &=& \gamma_{b}^6 - (\gamma_{b}^3)^2 - 6 \gamma_{b}^2
\gamma
_{b}^4 + 9(\gamma_{b}^2)^3 .
\end{eqnarray*}
\end{theorem}

As for second moments, these quantities denote the minimum variances,
which may be achieved for empirical estimators based on the true
$\varepsilon_{ij}$ and $b_i$, respectively.

%s2.4 ###
\vspace*{2pt}
\subsection{\texorpdfstring{Estimation of $\gamma_{\varepsilon}^4$ and $\gamma_b^4$}
{Estimation of gamma epsilon 4 and gamma b 4}}
\vspace*{2pt}
\label{sec24}

For $\gamma_{\varepsilon}^4$, we are also looking for a combination of
$f_m^4$'s such that the expectations include~$\gamma_{\varepsilon}^4$
but no other moments. First, from Lemma \ref{lem21}, we have
\begin{eqnarray*}
 \mathbb{E} f_4^4(i) &=& l_i \gamma_b^4 + 6 l_i \gamma_b^2 \gamma
_{\varepsilon}^2 + l_i \gamma_{\varepsilon}^4,\\
 \mathbb{E} f_3^4(i) &=& l_i^2 \gamma_b^4 + 3 l_i(l_i+1) \gamma_b^2
\gamma_{\varepsilon}^2 + l_i\gamma_{\varepsilon}^4,\\
 \mathbb{E} f_2^4(i) &=& l_i^3 \gamma_b^4 + (l_i^3 + 5l_i^2) \gamma_b^2
\gamma_{\varepsilon}^2 + l_i(l_i-1)(\gamma_{\varepsilon}^2)^2 +
l_i\gamma_{\varepsilon}^4
\end{eqnarray*}
and
\[
\mathbb{E} f_1^4(i) = l_i^4 \gamma_b^4 + 6 l_i^3
\gamma_b^2 \gamma_{\varepsilon}^2 + \mathbb{E} \biggl[ \biggl(\sum_j
\varepsilon_{ij}\biggr)^4\biggr] .
\]
Finally, we put
\[
f_5^4(i) = \Biggl[ \sum_{j=1}^{l_i} e_{ij}^2\Biggr]^2 .
\]
Clearly,
\begin{eqnarray}\label{eq19}
\mathbb{E} f_5^4(i) &=& \sum_{j=1}^{l_i} \sum_{k=1}^{l_i}
\mathbb{E} [e_{ij}^2 e_{ik}^2] = \sum_{j=1}^{l_i} \sum_{k=1}^{l_i}
\mathbb{E} [ (b_i + \varepsilon_{ij})^2 (b_i + \varepsilon_{ik})^2]
\nonumber\\ [-6pt]\\ [-6pt]
&=& l_i^2 \gamma_b^4 + l_i \gamma_{\varepsilon}^4 + (2 l_i^2 + 4 l_i)
\gamma_b^2 \gamma_{\varepsilon}^2 + (l_i^2 - l_i) (\gamma
_{\varepsilon
}^2)^2 .\nonumber
\end{eqnarray}
We now combine these expressions in a proper way. In particular, we
check that
\begin{eqnarray*}
&&\mathbb{E} \bigl[ (l_i^2- 2l_i +3) \bigl(l_i f_4^4(i) - 4 f_3^4(i)\bigr) + 6
l_i f_2^4(i) - 3f_1^4(i) - 3(2l_i - 3) f_5^4(i)\bigr] \\[2pt]
&&\quad = l_i(l_i - 1) (l_i - 2) (l_i - 3) \gamma_{\varepsilon}^4 .
\end{eqnarray*}
At first sight, the coefficients may look a little strange, but they
appear as solutions of linear equations incorporating
$\mathbb{E}f_1^4,\ldots, \mathbb{E}f_5^4$ such that all terms
involving moments other than $\gamma_{\varepsilon}^4$ vanish. Our
minimum variance estimator of $\gamma_{\varepsilon}^4$ thus becomes
\begin{eqnarray*}
\hat\gamma_{\varepsilon}^4 &=& N^{-1} \sum_{i=1}^n \frac
{1}{(l_i-1)(l_i-2)(l_i-3)}\\[2pt]
&&\hphantom{N^{-1} \sum_{i=1}^n}{}\times \Biggl\{ (l_i^2 - 2l_i + 3) \Biggl[ l_i \sum
_{j=1}^{l_i} \hat e_{ij}^4 - 4 \sum_{j=1}^{l_i} \hat e_{ij}^3 \sum
_{j=1}^{l_i} \hat e_{ij}\Biggr]\\[2pt]
&&\hspace*{-1pt}\hphantom{N^{-1} \sum_{i=1}^n\times \Biggl\{}{} + 6 l_i \Biggl( \sum_{j=1}^{l_i} \hat e_{ij}^2\Biggr) \Biggl( \sum_{j=1}^{l_i} \hat
e_{ij}\Biggr)^2 - 3\Biggl( \sum_{j=1}^{l_i} \hat e_{ij}\Biggr)^4 - 3 (2l_i - 3) \Biggl[
\sum_{j=1}^{l_i} \hat e_{ij}^2\Biggr]^2\Biggr\}.
\end{eqnarray*}
For $\gamma_b^4$, the relevant equation is
\[
\mathbb{E} [ f_1^4(i) - 6f_2^4(i) + 8f_3^4(i) - 6 f_4^4(i) + 3
f_5^4(i)] = l_i(l_i-1)(l_i-2)(l_i-3) \gamma_b^4,\
\]
giving us
\begin{eqnarray*}
&&\hat\gamma_b^4 = \frac{1}{n} \sum_{i=1}^n
\frac{1}{l_i(l_i-1)(l_i-2)(l_i-3)} \Biggl\{  \Biggl( \sum_{j=1}^{l_i} \hat
e_{ij}\Biggr)^4 - 6 \Biggl( \sum_{j=1}^{l_i} \hat e_{ij}^2\Biggr) \Biggl( \sum_{j=1}^{l_i}
\hat e_{ij}\Biggr)^2 \\[2pt]
&&\phantom{\hat\gamma_b^4 = \frac{1}{n} \sum_{i=1}^n
\frac{1}{l_i(l_i-1)(l_i-2)(l_i-3)} \Biggl\{}{} + 8 \Biggl( \sum_{j=1}^{l_i} \hat e_{ij}^3\Biggr) \Biggl( \sum_{j=1}^{l_i} \hat
e_{ij}\Biggr) - 6 \sum_{j=1}^{l_i} \hat e_{ij}^4 + 3 \Biggl( \sum_{j=1}^{l_i}
\hat
e_{ij}^2\Biggr)^2 \Biggr\}.
\end{eqnarray*}

\begin{theorem}\label{thm:6}
 Under the conditions of Theorem \ref{thm:1}, when
$\gamma_{\varepsilon}^8$ and $\gamma_b^8$ are finite, we have that
\[
N^{1/2} (\hat\gamma_{\varepsilon}^4 - \gamma_{\varepsilon}^4) \to
{\mathcal
N}_1(0, \mu_{\varepsilon}^4)
\]
and
\[
n^{1/2} (\hat\gamma_b^4 - \gamma_b^4) \to{\mathcal N}_1(0,\mu_b^4) ,
\]
where
\[
\mu_{\varepsilon}^4 = \gamma_{\varepsilon}^8 - (\gamma
_{\varepsilon
}^4)^2 - 8 \gamma_{\varepsilon}^3 \gamma_{\varepsilon}^5 + 16
\gamma
_{\varepsilon}^2 (\gamma_{\varepsilon}^3)^2
\]
and
\[
\mu_{b}^4 = \gamma_{b}^8 - (\gamma_{b}^4)^2 - 8 \gamma_{b}^3 \gamma
_{b}^5 + 16 \gamma_{b}^2 (\gamma_{b}^3)^2 .
\]
As in previous cases, $\mu_{\varepsilon}^4$ and $\mu_b^4$ are minimal
variances.
\end{theorem}

%s3 ###
\section{First step estimation}\label{sec3}

%s3.1 ###
\subsection{\texorpdfstring{Estimation of $\gamma_{\varepsilon}^3$ and $\gamma_b^3$}
{Estimation of gamma epsilon 3 and gamma b 3}}
\label{sec31}

In this section, we briefly discuss the fact that different choices of
estimating equations may lead to inefficiencies. These observations
eventually lead us to the efficient estimators discussed in the
previous section. For the third moments, recall that
\[
\mathbb{E} f_3^3(i) = l_i \gamma_b^3 + l_i \gamma_{\varepsilon}^3
\quad\mbox{and}\quad \mathbb{E} f_2^3(i) = l_i^2 \gamma_b^3 + l_i\gamma
_{\varepsilon}^3,
\]
from which
\[
l_i \gamma_{\varepsilon}^3 = \frac{1}{l_i - 1} [ l_i \mathbb{E}
f_3^3(i) - \mathbb{E} f_2^3(i)] .
\]
Summation over $1\le i\le n$ yields
\[
\gamma_{\varepsilon}^3 = \frac{\sum_{i=1}^n (1/(l_i - 1)) [ l_i
\mathbb{E} f_3^3(i) - \mathbb{E} f_2^3(i)]}{N} .
\]
If we replace the expectations by their sample analogs and the true
$e$'s by the residuals, then we come up with an estimator of
$\gamma_{\varepsilon}^3$ similar to that of Cox and Hall \cite{coxhall}, where
all~$l_i$'s are equal and there are no covariate effects:
\[
\hat\gamma_{\varepsilon}^{*3} = \frac{\sum_{i=1}^n (1/(l_i -1)) [
l_i \sum_{j=1}^{l_i} \hat e_{ij}^3 - ( \sum_{j=1}^{l_i} \hat e_{ij}^2)
( \sum_{j=1}^{l_i} \hat e_{ij}) ]}{N} .
\]
In the same way, we obtain
\[
\hat\gamma_b^{*3} = \frac{1}{n} \sum_{i=1}^n \frac{1}{l_i(l_i-1)}
\Biggl\{
\sum_{j=1}^{l_i} \hat e_{ij}^2 \sum_{j=1}^{l_i} \hat e_{ij} - \sum
_{j=1}^{l_i} \hat e_{ij}^3 \Biggr\} .
\]

To formulate limit results for $\hat\gamma_{\varepsilon}^{*3}$ and
$\hat\gamma_{b}^{*3}$, we recall that
\begin{eqnarray*}
\mu_{\varepsilon}^3 &=& \gamma_{\varepsilon}^6 - (\gamma
_{\varepsilon
}^3)^2 - 6 \gamma_{\varepsilon}^2 \gamma_{\varepsilon}^4 + 9(\gamma
_{\varepsilon}^2)^3,\\
\mu_{b}^3 &=& \gamma_{b}^6 - (\gamma_{b}^3)^2 - 6 \gamma_{b}^2
\gamma
_{b}^4 + 9(\gamma_{b}^2)^3
\end{eqnarray*}
and put
\[
\mu_{\varepsilon}^{*3} = \gamma_{\varepsilon}^6 - (\gamma
_{\varepsilon
}^3)^2 - 6 \gamma_{\varepsilon}^2 \gamma_{\varepsilon}^4 +
(4c+5)(\gamma
_{\varepsilon}^2)^3 + 4(\gamma_{\varepsilon}^2)^3 x_0'
\Sigma
^{-1} x_0 .
\]
Here,
\[
\bar{x}_n^* = N^{-1} \sum_{i=1}^n \sum_{j=1}^{l_i} x_{ij}
\]
and (as before)
\[
\bar{x}_n = \frac{1}{n} \sum_{i=1}^n \frac{1}{l_i} \sum_{j=1}^{l_i}
x_{ij} .
\]
The vector $x_0$ in $\mu_{\varepsilon}^{*3}$ equals
\[
x_0 = \lim_{n\to\infty} (\bar{x}_n^* - \bar{x}_n) ,
\]
while
\[
c = \lim_{n\to\infty} \frac{N}{n^2} \sum_{i=1}^n l_i^{-1},
\]
assuming that both limits exist.

A detailed qualitative interpretation of these quantities will be
deferred to the end of this section.

\begin{theorem}\label{thm:3}
 Under the conditions of Theorem \ref{thm:1}, when
$\gamma_{\varepsilon}^6$ and $\gamma_b^6$ are finite, we have that
%
%e13 ###
%
\begin{equation}\label{eq14}
 N^{1/2} (\hat\gamma_{\varepsilon}^{*3} -
\gamma_{\varepsilon}^3) \to{\mathcal N}_1\bigl(0, \mu_{\varepsilon}^{*3} +
4\gamma_b^2\bigl(\gamma_{\varepsilon}^4 -
(1-d)(\gamma_{\varepsilon}^2)^2\bigr)\bigr),
\end{equation}
where
\[
d = \lim_{n\to\infty} \biggl[ \frac{\sum_{i=1}^n l_i^2}{N} - \frac{\sum
_{i=1}^n l_i}{n} \biggr] .
\]
\end{theorem}

As to $\hat\gamma_b^{*3}$, we have that
%
%e14 ###
%
\begin{equation}\label{eq15}
 n^{1/2} [ \hat\gamma_b^{*3} - \gamma_b^3] \to{\mathcal N}_1
(0,\mu_b^3)\qquad \mbox{as } n \to\infty.
\end{equation}

\begin{rem}\label{rem31}
As in Section \ref{sec2}, the estimator in the $b$-case achieves the
minimum variance. It equals the variance of the moment estimator based
on the true but unknown $b_i$. In the $\varepsilon$-case, things are
less transparent. For example,\vadjust{\goodbreak} assume that $l_{ni} \equiv l_n^0$ are
all equal for $1\le i\le n$, a situation studied by Cox and Hall \cite{coxhall}. If
$l_n^0 \to\infty$, then $c = 1, x_0 = 0$ and $d = 0$. Hence,\vspace*{-1pt}
\[
\mu_{\varepsilon}^{\ast 3} =
\gamma_{\varepsilon}^6 -
(\gamma_{\varepsilon}^3)^2 - 6
\gamma_{\varepsilon}^2
\gamma_{\varepsilon}^4
+ 9(\gamma_{\varepsilon}^2)^3
= \mu_{\varepsilon}^3,\vspace*{-1pt}
\]
the variance of the (central) moment estimator based on the true
$\varepsilon_{ij}$. The total variance therefore becomes\vspace*{-1pt}
\[
\mu_{\varepsilon}^3 + 4 \gamma_b^3\bigl(\gamma_{\varepsilon}^4 -
(\gamma
_{\varepsilon}^2)^2\bigr),\vspace*{-1pt}
\]
which, by the Cauchy--Schwarz inequality, exceeds
$\mu_{\varepsilon}^3$. Hence, in this situation,
$\hat\gamma^{*3}_{\varepsilon}$ is inefficient.\vspace*{-1pt}
\end{rem}

\begin{rem}\label{rem32}
If $l_i = i^{a}$ with $0 < a < 1$, then $c = 1/(1-a^2) > 1$ becomes
large as $a \to1$. Hence, the quality of the Cox--Hall-type estimator
deteriorates in such situations. Worse than that, as our proofs reveal,
asymptotic normality may fail in situations where the limit $d$ is not
finite.\vspace*{-2pt}
\end{rem}

%s3.2 ###
\subsection{\texorpdfstring{Estimation of $\gamma_{\varepsilon}^4$ and $\gamma_b^4$}
{Estimation of gamma epsilon 4 and gamma b 4}}
\vspace*{-2pt}
\label{sec32}
%%2.4

For fourth moments, taking expectations of $f_4^4(i)$ and $f_3^4(i)$,
we again obtain\vspace*{-1pt}
\begin{equation}\label{eq16}
\mathbb{E} f_4^4(i)= l_i \gamma_b^4 + 6 l_i \gamma_b^2 \gamma
_{\varepsilon}^2 + l_i \gamma_{\varepsilon}^4\vspace*{-1pt}
\end{equation}
and\vspace*{-1pt}
\begin{equation} \label{eq17}
\mathbb{E} f_3^4(i)  = l_i^2 \gamma
_b^4 +
3 l_i (l_i + 1) \gamma_b^2 \gamma_{\varepsilon}^2 + l_i
\gamma_{\varepsilon}^4 ,\vspace*{-1pt}
\end{equation}
from which it follows that\vspace*{-1pt}
\[
\gamma_b^4 = \frac{\mathbb{E} f_3^4(i) - \mathbb
{E}f_4^4(i)}{l_i(l_i -
1)} - 3\gamma_b^2 \gamma_{\varepsilon}^2 .\vspace*{-1pt}
\]
Averaging over $1\le i\le n$ leads to the estimator\vspace*{-1pt}
\[
\hat{\hat{\gamma}}{}_b^{\ast4} = \frac{1}{n} \sum_{i=1}^n \frac
{1}{l_i(l_i - 1)} \Biggl\{ \sum_{j=1}^{l_i} \hat e_{ij}^3 \sum_{j=1}^{l_i}
\hat e_{ij} - \sum_{j=1}^{l_i} \hat e_{ij}^4 \Biggr\} - 3\hat\gamma_b^2
\hat\gamma_{\varepsilon}^2 ,\vspace*{-1pt}
\]
where $\hat\gamma_b^2$ and $\hat\gamma_{\varepsilon}^2$ were
studied in
Section \ref{sec22}. From (\ref{eq16}) and (\ref{eq17}) we immediately obtain\vspace*{-1pt}
\[
l_i \gamma_{\varepsilon}^4 = \frac{l_i \mathbb{E} f_4^4(i) -
\mathbb{E}
f_3^4(i)}{l_i - 1} - 3 l_i \gamma_b^2 \gamma_{\varepsilon}^2\vspace*{-1pt}
\]
and therefore to\vspace*{-1pt}
\[
\hat{\hat{\gamma}}{}_{\varepsilon}^{\ast4} = N^{-1} \sum_{i=1}^n
\frac
{1}{l_i - 1} \Biggl\{ l_i \sum_{j=1}^{l_i} \hat e_{ij}^4 - \sum_{j=1}^{l_i}
\hat e_{ij}^3 \sum_{j=1}^{l_i} \hat e_{ij}\Biggr\} - 3 \hat\gamma_b^2 \hat
\gamma_{\varepsilon}^2 .\vspace*{-1pt}
\]
Cox and Hall \cite{coxhall} also considered these estimators; however, we
have
discovered that the limit variances are\vadjust{\goodbreak} larger than those given in
their paper. Therefore, we propose the following modification. First,
recall
that
%
%e15 ###
%
\begin{equation}\label{eq18}
\mathbb{E}f_2^4(i) = l_i^3 \gamma_b^4 + (l_i^3 + 5l_i^2) \gamma_b^2
\gamma_{\varepsilon}^2 + l_i(l_i - 1) (\gamma_{\varepsilon}^2)^2 + l_i
\gamma_{\varepsilon}^4 .
\end{equation}
In addition to the $f_m^k(i)$ with $m \le k$, we again need
\[
f_5^4(i) = \Biggl[ \sum_{j=1}^{l_i} e_{ij}^2 \Biggr]^2 .
\]
It follows from (\ref{eq18}) and (\ref{eq19}) that
%
%e16 ###
%
\begin{equation}\label{eq20}
\mathbb{E} [ f_2^4(i) - f_5^4(i)] = (l_i^3 - l_i^2)
\gamma
_b^4 + (l_i^3 + 3l_i^2 - 4l_i) \gamma_{\varepsilon}^2 \gamma_b^2 .
\end{equation}
To estimate $\gamma_{\varepsilon}^4$, we are looking for a linear
combination of (\ref{eq16}), (\ref{eq17}) and (\ref{eq20}) so that the
terms~$\gamma_b^4$ and $\gamma_{\varepsilon}^2\gamma_b^2$ cancel out.
In fact, it is easily seen that
\[
\mathbb{E} [ (2 l_i^2 - l_i) f_4^4(i) - (5l_i - 4) f_3^4(i) + 3
f_2^4(i) - 3f_5^4(i)] = 2 l_i (l_i - 1)(l_i - 2)\gamma_{\varepsilon}^4.
\]
The corresponding estimator of $\gamma_{\varepsilon}^4$ becomes
\begin{eqnarray*}
&&\hat\gamma_{\varepsilon}^{\ast4} = N^{-1} \sum_{i=1}^n \frac{1}{2(l_i
- 1)(l_i - 2)} \Biggl\{ (2 l_i^2 -l_i) \sum_{j=1}^{l_i} \hat e_{ij}^4
- (5l_i - 4) \sum_{j=1}^{l_i} \hat e_{ij}^3 \sum_{j=1}^{l_i} \hat
e_{ij} \\
&&\hphantom{\hat\gamma_{\varepsilon}^{\ast4} = N^{-1} \sum_{i=1}^n \frac{1}{2(l_i
- 1)(l_i - 2)} \Biggl\{}{} + 3 \Biggl( \sum_{j=1}^{l_i} \hat e_{ij}^2\Biggr) \Biggl( \sum_{j=1}^{l_i} \hat
e_{ij}\Biggr)^2 - 3\Biggl[ \sum_{j=1}^{l_i} \hat e_{ij}^2\Biggr]^2 \Biggr\} .
\end{eqnarray*}
Following this idea, we also get an estimator of $\gamma_b^4$.
Subtracting (\ref{eq16}) from (\ref{eq17}), we obtain
\[
\mathbb{E} [f_3^4(i) - f_4^4(i)] = (l_i^2 - l_i) \gamma_b^4 + 3(l_i^2 -
l_i) \gamma_b^2 \gamma_{\varepsilon}^2 .
\]
Together with (\ref{eq20}), this yields
\[
3\mathbb{E} [f_2^4(i) - f_5^4(i)] - (l_i + 4) \mathbb{E} [f_3^4(i) -
f_4^4(i)] = 2 l_i (l_i - 1) (l_i - 2) \gamma_b^4
\]
and therefore
\begin{eqnarray*}
&&\hat\gamma_b^{\ast4} = \frac{1}{n} \sum_{i=1}^n \frac{1}{2
l_i(l_i-1)(l_i-2)} \Biggl\{ 3\Biggl( \sum_{j=1}^{l_i} \hat e_{ij}^2\Biggr) \Biggl(
\sum
_{j=1}^{l_i} \hat e_{ij}\Biggr)^2 - 3\Biggl( \sum_{j=1}^{l_i} \hat
e_{ij}^2\Biggr)^2 \\
&&\hphantom{\hat\gamma_b^{\ast4} = \frac{1}{n} \sum_{i=1}^n \frac{1}{2
l_i(l_i-1)(l_i-2)} \Biggl\{}{} - (l_i + 4) \Biggl( \sum_{j=1}^{l_i} \hat e_{ij}^3\Biggr) \Biggl( \sum_{j=1}^{l_i}
\hat e_{ij}\Biggr) + (l_i + 4) \sum_{j=1}^{l_i} \hat e_{ij}^4 \Biggr\} .
\end{eqnarray*}
In the following theorem, we summarize the main results on the limit
distributions of~$\hat\gamma_{\varepsilon}^{\ast4}$
and~$\hat\gamma_b^{\ast4}$.
\begin{theorem}\label{thm:4}
Under the conditions of Theorem \ref{thm:3}, when
$\gamma_{\varepsilon}^8$ and $\gamma_b^8$ are finite, we have that
\[
N^{1/2} (\hat\gamma_{\varepsilon}^{\ast4} - \gamma_{\varepsilon}^4)
\to{\mathcal N}_1\bigl(0, \mu_{\varepsilon}^{\ast4} + \tfrac{9}{4} \gamma_b^2
[\gamma_{\varepsilon}^6 - (1-d) (\gamma_{\varepsilon}^3)^2 -
6\gamma_{\varepsilon}^2 \gamma_{\varepsilon}^4 +
9(\gamma_{\varepsilon}^2)^3]\bigr)
\]
and
\[
n^{1/2} [\hat\gamma_b^{\ast4} - \gamma_b^4] \to{\mathcal N}_1 (0, \mu_b^4)\qquad
\mbox{as } n \to\infty,
\]
where, again,
\[
\mu_{b}^4 = \gamma_{b}^8 - (\gamma_{b}^4)^2 - 8 \gamma_{b}^3 \gamma
_{b}^5 + 16 \gamma_{b}^2 (\gamma_{b}^3)^2
\]
and
\[
\mu_{\varepsilon}^{\ast4} = \gamma_{\varepsilon}^8 - (\gamma
_{\varepsilon}^4)^2 - 8 \gamma_{\varepsilon}^3 \gamma_{\varepsilon
}^5 +
\bigl( \tfrac{9}{4}c + \tfrac{55}{4}\bigr) \gamma_{\varepsilon}^2
(\gamma_{\varepsilon}^3)^2 + \tfrac{9}{4} \gamma_{\varepsilon}^2
(\gamma
_{\varepsilon}^3)^2 x'_0 \Sigma^{-1} x_0 .
\]
\end{theorem}

\begin{rem}
Our earlier Remarks \ref{rem31} and \ref{rem32} also apply to fourth moments. This more
or less led us to look for the new estimators studied in Section \ref{sec2}.
\end{rem}

\begin{rem}
We will now discuss the results of this paper in a qualitative way.
Suppose that all the $b_i$'s and $\varepsilon_{ij}$'s are known to the
observer. Then, rather than computing residuals, they could be used
directly to nonparametrically estimate the (central) moments of
interest. Simple computations then show that the variances of these
estimators equal $\mu_b^k$ and $\mu_{\varepsilon}^k$, respectively. In
the case where only residuals are available, the improper weighting in
$\hat\gamma_{\varepsilon}^{\ast k}$ yields variances which heavily
depend on the design (via $x'_0 \Sigma^{-1}x_0)$, the group sizes (via
the constants $c$ and~$d$) and the noise variables $b_i$ (via
$\gamma_b^2$). In such a situation, Theorems \ref{thm:5} and \ref{thm:6} provide new
estimators which also attain the minimum variance in the
$\varepsilon$-case and are less vulnerable to the model design.
\end{rem}

%s4 ###
\section{Simulation study}\label{sec4}
To demonstrate the usefulness of our estimation procedures, a small
simulation study will be carried out. The data sets are generated from
the model (\ref{eq2.1}) with $\alpha=1$ and $\beta=(1,2)'$. To estimate the
model parameters and the third and fourth moments using the methods
developed in this paper, the group values are randomly drawn from a
Poisson distribution with mean $5$. The design matrices are generated
from a zero-mean normal distribution with covariance matrix $ \left({
1 \enskip 0.8}\atop{0.8 \enskip 1}\right)$. For the random effects $b_i$ and the errors $\varepsilon_{ij}$, we
consider the following five cases:
\begin{longlist}
\item[(a)]
$\varepsilon_{ij}\sim_{\mathrm{i.i.d.}}0.5{\mathcal N}_1(0,1)$ and $b_i \sim
_{\mathrm{i.i.d.}}0.5{\mathcal N}_1(0,1)$;
\item[(b)]
$\varepsilon_{ij}\sim_{\mathrm{i.i.d.}}0.5{\mathcal N}_1(0,1)$ and $b_i \sim
_{\mathrm{i.i.d.}}0.5t(8)$;
\item[(c)]
$\varepsilon_{ij}\sim_{\mathrm{i.i.d.}}0.5{\mathcal N}_1(0,1)$ and $b_i \sim
_{\mathrm{i.i.d.}}0.5\Gamma(1,1)-0.5$;
\item[(d)]
$\varepsilon_{ij}\sim_{\mathrm{i.i.d.}}0.5t(8)$ and $b_i \sim_{\mathrm{i.i.d.}}0.5t(8)$;
\item[(e)]
$\varepsilon_{ij}\sim_{\mathrm{i.i.d.}}0.5t(8)$ and $b_i \sim_{\mathrm{i.i.d.}}0.5\Gamma
(1,1)-0.5$.
\end{longlist}
The true values of the $2$nd--$4$th moments of the errors and random
effects are given in Table \ref{tab1}. ${\mathcal N}_1$, $\Gamma$ and $t$
correspond to the normal, gamma and $t$ distributions,
respectively.\vadjust{\goodbreak}

\begin{table}
\tablewidth=220pt
\caption{The true values of the $2$nd--$4$th moments of the random and group
errors}\label{tab1}
\vspace*{-3pt}
\begin{tabular*}{220pt}{@{\extracolsep{\fill}}llll@{}}
\hline
 c.d.f                & $2$nd & $3$rd & $4$th\\
 \hline
 $0.5N(0,1)$          & 0.25  & 0     & 0.1875 \\
 $0.5t(8)$            & 0.333 & 0     &0.5 \\
 $0.5\Gamma(1,1)-0.5$ & 0.25  &0.25   &0.5625\\
 \hline
\end{tabular*}
\vspace*{-6pt}
\end{table}

The following simulation results are based on $1000$ samples of data
$\{(x_{ij},y_{ij})\dvtx i=1,\ldots,n, j=1,\ldots,l_i\}$ with $n=50,100$.
The estimated mean, standard deviation and root mean squared error of
the estimators suggested above are reported in Tables \ref{tab2} and \ref{tab3}. Table
\ref{tab2}
presents the results for the model parameters and second moments. For
the purposes of comparison, we also include the results for the MLE.
Table \ref{tab3} presents the results for the minimum variance estimators of
the third and fourth moments.

\begin{table}
\tabcolsep=0pt
 \caption{The results for $\hat\alpha$, $\hat\beta$, $\hat\gamma_{\varepsilon}^2$ and $\hat\gamma_{b}^2$ in cases (a)--(e) (the numbers in
          brackets correspond to the MLE)}\label{tab2}
\begin{tabular*}{\textwidth}{@{\extracolsep{\fill}}llllllll@{}}
\hline
Case&$n$&Result&$\hat{\alpha}$&$\hat\beta_1$&$\hat\beta_2$&$\hat\gamma_{\varepsilon}^2$&$\hat\gamma_{b}^2$ \\
\hline
(a)&\phantom{0}$50$&mean&$1.0004\,(1.0005)$&$ 0.9986\,(0.9984)$&$2.0009\,(2.0012)$ &$0.2497\,(0.2480)$  &$0.2451\,(0.2449)$  \\
&&std &$0.0748\,(0.0749)$&$ 0.0428\,(0.0424)$ &$0.0430\,( 0.0427)$ &$0.0177\,(0.0177)$ &$0.0550\,(0.0549)$\\
 &&rmse &$0.0748\,(0.0749)$& $ 0.0428\,(0.0424)$&$0.0430\,(
0.0427)$ &$0.0178\,(0.0178)$&$0.0552\,(0.0552)$\\
&$100$&mean&$0.9965\,(0.9964)$ &$1.0000\,(0.9999)$&$1.9995\,(1.9996)$& $0.2496\,(0.2496)$&$0.2476\,(0.2473)$  \\
 && std&$0.0544\,(0.0544)$& $ 0.0307\,(0.0304)$&$0.0302\,( 0.0300)$&$0.0128\,(0.0128)$& $0.0392\,(0.0393)$\\
  &&rmse&$0.0545\,(0.0546)$& $ 0.0307\,(0.0304)$&$0.0302\,( 0.0300)$&$0.0128\,(0.0128)$&$0.0393\,(0.0394)$
  \\ [3pt]
(b)&\phantom{0}$50$&mean&$1.0022\,(1.0023)$ &$ 1.0004\,(1.0002)$ &$2.0010\,(2.0014)$ &$0.2490\,(0.2490)$ &$0.3301\,(0.3296)$   \\
 &&std&$0.0860\,(0.0860)$& $ 0.0438\,(0.0438)$ &$0.0427\,( 0.0427)$&$0.0184\,(0.0184)$ &$0.0978\,(0.0978)$\\
 && rmse&$0.0861\,(0.0860)$& $ 0.0438\,(0.0438)$ &$0.0427\,( 0.0427)$&$0.0185\,(0.0184)$&$0.0979\,(0.0979)$ \\
&$100$&mean&$1.0026\,(1.0026)$ &$0.9991\,(0.9991)$ &$2.0005\,(2.0004)$ &$0.2489\,(0.2489)$&$ 0.3322\,(0.3319 )$  \\
 && std &$0.0598\,(0.0598)$& $ 0.0299\,(0.0296)$ &$0.0299\,( 0.0296)$&$0.0128\,(0.0127)$&$ 0.0652\,(0.0651)$ \\
 && rmse&$0.0598\,(0.0598)$& $ 0.0299\,(0.0296)$ &$0.0299\,( 0.0296)$&$0.0128\,(0.0128)$&$0.0652\,(0.0651)$
 \\ [3pt]
 (c)&\phantom{0}$50$&mean&$0.9988\,(0.9988)$&$ 1.0004\,(1.0003)$ &$1.9962\,(1.9960)$ &$0.2485\,(0.2485)$&$0.2421\,(0.2418)$  \\
 &&std &$0.0754\,(0.0754)$& $ 0.0430\,(0.0426)$ &$0.0418\,( 0.0415)$&$0.0178\,(0.0178)$ &$0.0981\,(0.0981)$\\
 &&rmse &$0.0754\,(0.0754)$&$ 0.0430\,(0.0426)$ &$0.0420\,( 0.0417)$&$0.0179\,(0.0179)$&$0.0984\,(0.0984)$ \\
&$100$&mean&$0.9974\,(0.9974)$ &$0.9993\,(0.9993)$ &$2.0005\,(2.0007)$ &$0.2481\,(0.2481)$&$ 0.2437\,(0.2440)$  \\
&&std &$0.0527\,(0.0528)$& $ 0.0298\,(0.0294)$ &$0.0302\,( 0.0298)$&$0.0121\,(0.0121)$ &$0.0731\,(0.0730)$\\
&&rmse&$0.0527\,(0.0528)$&$ 0.0298\,(0.0294)$&$0.0302\,( 0.0298)$
&$0.0123\,(0.0123)$&$0.0733\,(0.0733)$\\ [3pt]
 (d)&\phantom{0}$50$&mean&$1.0008\,(1.0007)$&$ 0.9992\,(0.9992)$ &$2.0001\,(2.0011)$&$0.3298\,(0.3297)$&$0.3239 \,( 0.3232)$    \\
 &&std &$0.0843\,(0.0842)$&$ 0.0479\,(0.0474)$ &$0.0484\,( 0.0479)$&$0.0294\,(0.0294)$&$0.0896\,(0.0894)$ \\
&&rmse&$0.0843\,(0.0842)$&
$0.0479\,(0.0474)$&$0.0484\,(0.0480)$&$0.0296\,(0.0296)$&$0.0900\,(0.0900)$\\
  &$100$&mean&$1.0004\,(1.0003)$ &$1.0006\,(1.0004)$ &$1.9996\,(1.9996)$ &$0.3315\,(0.3315)$&$ 0.3305\,(0.3303 )$   \\
  &&std &$0.0618\,(0.0618)$& $ 0.0336\,(0.0334)$ &$0.0351\,( 0.0350)$ &$0.0213\,(0.0213)$&$ 0.0697\,(0.0697)$ \\
  &&rmse&$0.0618\,(0.0618)$& $ 0.0336\,(0.0334)$ &$0.0351\,( 0.0350)$&$0.0214\,(0.0214)$&$0.0697\,(0.0697)$
  \\ [3pt]
(e)&\phantom{0}$50$&mean&$1.0025\,(1.0023)$ &$ 1.0016\,(1.0012)$ &$1.9977\,(1.9977)$&$0.3311\,(0.3313)$&$ 0.2433\,(0.2425)$   \\
&&std &$0.0670\,(0.0670)$& $ 0.0488\,(0.0488)$ &$0.0494\,( 0.0492)$&$0.0297\,(0.0298)$&$0.1010\,(0.1014)$ \\
&&rmse&$0.0670\,(0.0671)$& $ 0.0489\,(0.0488)$
&$0.0495\,(0.0493)$&$0.0297\,(0.0298)$&$0.1012\,(0.1017)$
\\
&$100$&mean&$1.0007\,(1.0006)$ &$0.9994\,(0.9993)$ &$2.0011\,(2.0010)$&$0.3323\,(0.3323)$&$ 0.2446\,(0.2442 )$   \\
&& std &$0.0537\,(0.0536)$& $ 0.0348\,(0.0344)$&$0.0335\,( 0.0330)$ &$0.0220\,(0.0220)$&$ 0.0708\,(0.0707)$ \\
&& rmse&$0.0537\,(0.0536)$&$0.0348\,(0.0344)$&$0.0335\,(
0.0330)$&$0.0220\,(0.0220)$&$0.0710\,(0.0709)$ \\
\hline
\end{tabular*}
\end{table}

In Table \ref{tab2}, the comparison with the MLE shows that our estimators are
very competitive, although such a comparison is actually in favor of
the MLE when we assume that the distribution is parametric. In fact,
empirical studies in the literature also show that the assumption
concerning the distribution of the random effects hardly influences
the parameter estimates; see Butler and Louis \cite{butlou} and Verbeke and Lesaffre \cite{verles} for
details. This indicates that the estimation of moments for the model
parameters performs very well.\looseness=-1

\begin{table}
\caption{The results for $\hat\gamma_{\varepsilon}^k$ and $\hat\gamma_{b}^k$ $(k=3,4)$ in cases
(a)--(e)}\label{tab3}
\begin{tabular*}{\textwidth}{@{\extracolsep{\fill}}llld{2.4}ld{2.4}l@{}}
\hline
\multicolumn{1}{@{}l}{Case}&\multicolumn{1}{l}{$n$}&\multicolumn{1}{l}{Result}
&\multicolumn{1}{l}{$\hat\gamma_{\varepsilon}^3$}&
\multicolumn{1}{l}{$\hat\gamma_{\varepsilon}^4$}&\multicolumn{1}{l}{$\hat\gamma_{b}^3$}
&\multicolumn{1}{l@{}}{$\hat\gamma_{b}^4$}\\
\hline
(a)&\phantom{0}50&mean&  -0.0004   & 0.1852&  -0.0003  & 0.1813\\
&&std &  0.0170 &0.0325&  0.0519 &0.0974 \\
&&rmse& 0.0170 &0.0326& 0.0519&0.0976\\
 &100&mean&-0.0002  & 0.1867&  0.0005 & 0.1835 \\
&& std&  0.0125 &0.0234&  0.0362&0.0688\\
&& rmse&  0.0125 &0.0234&  0.0362 &0.0689\\ [3pt]
 (b)&\phantom{0}50&mean&  0.0002  & 0.1866 & -0.0003  & 0.5054   \\
&&std&  0.0171 &0.0331&  0.1948 &0.8577  \\
&&rmse& 0.0171 &0.0331&
0.1948 &0.8578 \\
&100&mean & 0.0002  & 0.1858 & 0.0076  & 0.5003   \\
&&std&  0.0121 & 0.0239 &  0.1276&0.4835 \\
&&rmse& 0.0121 & 0.0240 &  0.1278  & 0.4835 \\ [3pt]
(c)&\phantom{0}50&mean&  0.0007  & 0.1864  &  0.2290  & 0.4962  \\
&&std &  0.0174 & 0.0335 &  0.2235 &0.7649 \\
&& rmse&  0.0174 &0.0335&
0.2244 & 0.7678 \\
&100&mean&  0.0002  & 0.1848 &  0.2380  & 0.5317   \\
&& std &  0.0123 &0.0216&  0.1778 &0.6445 \\
&& rmse& 0.0123&0.0216&  0.1782 &0.6445  \\ [3pt]
(d)&\phantom{0}50&mean& 0.0001   & 0.4754& 0.0033   & 0.4477    \\
&&std&  0.0557 &0.1832&  0.1520&0.4510  \\
&&rmse &  0.0557 &0.1848&
0.1520 & 0.4542 \\
&100&mean& 0.0007   & 0.4862&  0.0014  & 0.4796   \\
&& std&  0.0403 &0.1611&  0.1098 &0.4100 \\
&&rmse& 0.0403&0.1617&  0.1098 &0.4102 \\ [3pt]
(e)&\phantom{0}50&mean&       -0.0011 & 0.4832&  0.2278  & 0.5007   \\
&&std&    0.0578 &0.2068 &  0.2311 &0.8052 \\
&&rmse& 0.0578 &0.2075&
0.2322 & 0.8075 \\
&100&mean&  0.0001   & 0.4979 &  0.2385  & 0.5355   \\
&&std& 0.0427&0.1726&  0.1753 &0.6569 \\
&&rmse& 0.0427 &0.1726 & 0.1757 &0.6575\\
\hline
\end{tabular*}
\vspace*{-3pt}
\end{table}

\begin{appendix}\label{appendix}
%s5 ###
\vspace*{-5pt}\section*{Appendix}\vspace*{-5pt}
\setcounter{section}{5}
\renewcommand{\theequation}{\arabic{section}.\arabic{equation}}

\begin{pf*}{Proof of Theorem \ref{thm:1}} We first study $\hat\beta$. It follows
from (\ref{eq2.1}) and (\ref{eq2.3}) that\vspace*{-2pt}
%
%e17 ###
%
\begin{equation}\label{eq5.1}
\hat\beta- \beta= \frac{1}{\sum_{i=1}^n l_i} \sum_{i=1}^n \sum
_{j=1}^{l_i} \hat\Sigma_n^{-1} (x_{ij} - \bar{x}_{i\cdot})
\varepsilon_{ij} ,\vspace*{-2pt}
\end{equation}
where $\hat{\Sigma}_n$ is given in (\ref{eq2.5}). To show
(\ref{eq2.9}), we fix $a \in\mathbb{R}^p$. It suffices to prove that\vspace*{-2pt}
\[
N^{1/2} a'(\hat\beta- \beta) \to{\mathcal N}_1(0,\gamma_{\varepsilon}^2
a'\Sigma^{-1} a)\qquad \mbox{in distribution}.\vspace*{-2pt}
\]
Since, according to (\ref{eq5.1}), $\hat\beta- \beta$ is a sum of
zero-mean independent random vectors, it remains to check the variance
and verify Lindeberg's condition. The variance of $N^{1/2} a'(\hat
\beta
- \beta)$ equals\vspace*{-2pt}
\[
\frac{\gamma_{\varepsilon}^2}{N} \sum_{i=1}^n \sum_{j=1}^{l_i} [
a'\hat{\Sigma}_n^{-1} (x_{ij} - \bar{x}_{i\cdot})]^2 =
\gamma_{\varepsilon}^2 a'\hat{\Sigma}_n^{-1} a \to\gamma
_{\varepsilon}^2 a'\Sigma^{-1} a ,\vspace*{-2pt}
\]
by (\ref{eq2.6}). To verify Lindeberg's condition, we first fix
$\delta
> 0$. The Lindeberg function then equals\vspace*{-2pt}
\[
L_n(\delta) = N^{-1} \sum_{i=1}^n \sum_{j=1}^{l_i} [ a'
\hat{\Sigma}_n^{-1} (x_{ij} - \bar{x}_{i\cdot})]^2 \int_{\{
|a'\hat\Sigma_n^{-1} (x_{ij} - \bar{x}_{i\cdot})| |\varepsilon_{ij}|
\ge\delta\sqrt{\Sigma l_i}\}} \varepsilon_{ij}^2\,\mathrm{d}\mathbb{P}.\vspace*{-2pt}
\]
Recall that, by (\ref{eq2.8}),
\[
C_n = \frac{\max_{i,j} \| x_{ij} - \bar{x}_{i\cdot}\|}{N^{1/2}} \to0\qquad
\mbox{as } n \to\infty.
\]
We conclude, by the Cauchy--Schwarz inequality and (\ref{eq2.6}),
that
\[
L_n(\delta) \le a' \hat{\Sigma}_n^{-1} a
\int_{\{\| a'\hat\Sigma_n^{-1}\| |\varepsilon| \ge\delta
C_n^{-1}\}} \varepsilon^2\,\mathrm{d}\mathbb{P} \to0,
\]
as required. This proves (\ref{eq2.9}). For $\hat\alpha$, we
immediately get from (\ref{eq2.4}) that
\[
\hat\alpha- \alpha= \frac{1}{n} \sum_{i=1}^n b_i + \bar
{\varepsilon}
- \frac{1}{n} \sum_{i=1}^n \bar{x}_{i\cdot}' (\hat\beta- \beta) .
\]
From (\ref{eq2.7}), it follows that $n^{1/2} \bar{\varepsilon} \to0$
in squared mean and hence in probability.

Furthermore, by (\ref{eq2.7})--(\ref{eq2.9}),
\[
n^{-1/2} \sum_{i=1}^n \bar{x}'_{i\cdot} (\hat\beta- \beta) \to0\qquad
\mbox{in probability}.\vadjust{\goodbreak}
\]
Summarizing,
\[
n^{1/2} (\hat\alpha- \alpha) = n^{-1/2} \sum_{i=1}^n b_i +
\mathrm{o}_{\mathbb
{P}}(1) \to{\mathcal N}_1(0, \gamma_b^2) .
\]
This shows (\ref{eq2.10}) and thereby completes the proof of Theorem \ref{thm:1}.
\end{pf*}

\begin{pf*}{Proof of Theorem \ref{thm:2}} From the definition of
$\hat\gamma_{\varepsilon}^2$ and (\ref{eq2.11}), we readily get
\begin{eqnarray*}
N \hat\gamma_{\varepsilon}^2 &=& \sum_{i=1}^n \frac{1}{l_i - 1}
\Biggl\{
l_i \sum_{j=1}^{l_i} \bigl[ \varepsilon_{ij}^2 + \bigl(z_{ij}'(\beta-
\hat\beta)\bigr)^2 + 2\varepsilon_{ij}
z'_{ij} (\beta-\hat\beta)\bigr] \\
&&\hphantom{\sum_{i=1}^n \frac{1}{l_i - 1}
\Biggl\{}{}- \sum_{j=1}^{l_i} \sum_{k=1}^{l_i} [ \varepsilon
_{ij} + z_{ij}'(\beta- \hat\beta)] [ \varepsilon_{ik} +
z_{ik}'(\beta- \hat\beta)] \Biggr\}\\
&=& \sum_{i=1}^n \sum_{j=1}^{l_i} \varepsilon_{ij}^2 - \sum_{i=1}^n
\frac{1}{l_i-1} \sum_{j\not= k} \varepsilon_{ij} \varepsilon_{ik}
+ 2
\sum_{i=1}^n \sum_{j=1}^{l_i} \varepsilon_{ij} z'_{ij}(\beta- \hat
\beta) \\
&&{} + \sum_{i=1}^n \frac{l_i}{l_i-1} \sum_{j=1}^{l_i} \bigl(z'_{ij}(\beta-
\hat\beta)\bigr)^2 - \sum_{i=1}^n \frac{1}{l_i-1} \sum_{j=1}^{l_i} \sum
_{k=1}^{l_i} [ z'_{ij}(\beta- \hat\beta)] [
z'_{ik}(\beta- \hat\beta)]\\
&&{} - \sum_{i=1}^n \frac{2}{l_i - 1} \sum_{j\not= k} \varepsilon_{ij}
z'_{ik} (\beta- \hat\beta)\\
& \equiv& I - \mathit{II} + \mathit{III} + \mathit{IV} - V - \mathit{VI}.
\end{eqnarray*}
Of these six terms, only the first will be a leading term, while the
others are remainders. For example, II is a sum of centered independent
random variables with variance
\[
2(\gamma_{\varepsilon}^2)^2 \sum_{i=1}^n \frac{l_i}{l_i-1} \le4n
(\gamma_{\varepsilon}^2)^2 .
\]
We conclude, in view of (\ref{eq2.7}), that
\[
N^{-1/2} \mathit{II} \to0\qquad \mbox{in probability}.
\]
To show the same for III, it suffices to prove, because of
(\ref{eq2.9}), that
\[
N^{-1} \sum_{i=1}^n \sum_{j=1}^{l_i} \varepsilon_{ij} z_{ij} \to0\qquad
\mbox{in probability}.
\]
Again, this is a sum of centered random vectors with covariance
\[
\gamma_{\varepsilon}^2 \biggl[ \frac{\sum_{i=1}^n l_i(\bar{x}_{i\cdot} -
\bar{x}) (\bar{x}_{i\cdot} - \bar{x})'}{N^2} + \frac{
\hat\Sigma_n}{N}\biggr] \to0 .
\]
Similarly, the convergence of $N^{-1/2}\mathit{IV}$, $N^{-1/2}V$ and
$N^{-1/2}\mathit{VI}$ to zero follows from~(\ref{eq2.6})--(\ref{eq2.9}). All
together, this shows that\vspace*{-3pt}
\[
N^{1/2} [ \hat\gamma_{\varepsilon}^2 - \gamma_{\varepsilon}^2 ] =
N^{-1/2} \Biggl[ \sum_{i=1}^n \sum_{j=1}^{l_i} (\varepsilon_{ij}^2 -
\gamma
_{\varepsilon}^2)\Biggr] + \mathrm{o}_{\mathbb{P}}(1)\vspace*{-3pt}
\]
and hence (\ref{eq2.12}), by a simple application of the central limit
theorem. To show (\ref{eq2.13}), we note that\vspace*{-3pt}
\[
\hat\gamma_b^2 = \frac{1}{n} \sum_{i=1}^n \frac{1}{l_i(l_i-1)}
\sum
_{j\not= k} \hat e_{ij}\hat e_{ik} .\vspace*{-3pt}
\]
Again using (\ref{eq2.11}) and applying similar arguments to those
used before, we obtain\vspace*{-3pt}
\[
\hat\gamma_b^2 = \frac{1}{n} \sum_{i=1}^n b_i^2 + \mathrm{o}_{\mathbb
{P}}(n^{-1/2}) ,\vspace*{-3pt}
\]
from which it follows that (\ref{eq2.13}) holds.\vspace*{-3pt}
\end{pf*}

\begin{pf*}{Proof of Theorem \ref{thm:5}} We first deal with $\hat\gamma_b^3$.
Simple algebraic manipulations yield\vspace*{-3pt}
\[
\hat\gamma_b^3 = \frac{1}{n} \sum_{i=1}^n \frac{1}{l_i(l_i-1)(l_i-2)}
\sum_{j\not=k\not=l} \hat e_{ij} \hat e_{ik} \hat e_{il} .\vspace*{-3pt}
\]
Expanding $\hat e_{ij}$ into $\hat e_{ij} = \varepsilon_{ij} + b_i -
(\bar{b} + \bar{\varepsilon}) + z'_{ij}(\beta-\hat\beta)$, we may
again neglect all contributions involving the $z'_{ij}(\beta-
\hat\beta)$. Hence, up to an $\mathrm{o}_{\mathbb{P}}(n^{-1/2})$ term,\vspace*{-3pt}
\begin{eqnarray*}
\hat\gamma_b^3 &=& \frac{1}{n} \sum_{i=1}^n \frac{1}{l_i(l_i-1)(l_i-2)}
\sum_{j\not= k \not= l} \varepsilon_{ij} \varepsilon
_{ik}\varepsilon
_{il} \\[-3pt]
&&{} + \frac{1}{n} \sum_{i=1}^n \frac{1}{l_i(l_i-1)} \sum_{j\not= k}
[3b_i \varepsilon_{ij} \varepsilon_{ik} - 3(\bar{b} + \bar
{\varepsilon
}) \varepsilon_{ij} \varepsilon_{ik}] \\[-3pt]
&&{} + \frac{1}{n} \sum_{i=1}^n \frac{1}{l_i} \sum_{j} [3b_i^2
\varepsilon
_{ij} - 6b_i(\bar{b} + \bar{\varepsilon}) \varepsilon_{ij} + 3(\bar{b}
+ \bar{\varepsilon})^2 \varepsilon_{ij}] \\[-3pt]
&&{} + \frac{1}{n} \sum_{i=1}^n [b_i^3 - 3(\bar{b} + \bar{\varepsilon})
b_i^2] + 3 \bar{b} (\bar{b} + \bar{\varepsilon})^2 - (\bar{b} +
\bar
{\varepsilon})^3 .\vspace*{-3pt}
\end{eqnarray*}
Under the assumptions of the theorem, the first three sums are
negligible, as are the last two terms. Hence,\vspace*{-3pt}
\begin{eqnarray*}
n^{1/2} [\hat\gamma_b^3 - \gamma_b^3]
&=& n^{-1/2} \sum_{i=1}^n [b_i^3 - \gamma_b^3 - 3(\bar{b} + \bar
{\varepsilon}) b_i^2] + \mathrm{o}_{\mathbb{P}}(1)\\[-3pt]
&=& n^{-1/2} \sum_{i=1}^n [b_i^3 - \gamma_b^3 - 3\gamma_b^3 b_i] +
\mathrm{o}_{\mathbb{P}}(1) .\vspace*{-3pt}\vadjust{\goodbreak}
\end{eqnarray*}
The conclusion for $\hat\gamma_b^3$ now readily follows from the
central limit theorem. For $\hat\gamma_{\varepsilon}^3$, we may write
\begin{eqnarray*}
N \hat\gamma_{\varepsilon}^3
&=& \sum_i \frac{2}{(l_i - 1)(l_i - 2)} \sum_{j\not= k \not= l}
\hat
e_{ij} \hat e_{ik} \hat e_{il} \\
&&{} - \sum_i \frac{3}{l_i - 1} \sum_{j \not= k} \hat e_{ij}^2 \hat
e_{ik} + \sum_{i} \sum_{j} \hat e_{ij}^3 .
\end{eqnarray*}
If we once again ignore the higher order terms of $z'_{ij}(\beta-
\hat\beta)$, we find that in the expansion of
$\hat\gamma_{\varepsilon}^3$, we have
\begin{eqnarray*}
N \hat\gamma_{\varepsilon}^3
&=& \sum_i \frac{2}{(l_i - 1)(l_i - 2)} \sum_{j\not= k \not= l}
\varepsilon_{ij} \varepsilon_{ik} \varepsilon_{il} \\
&&{} - \sum_i \frac{3}{l_i - 1} \sum_{j \not= k} \varepsilon_{ij}^2
\varepsilon_{ik} + \sum_{i} \sum_{j} \varepsilon_{ij}^3 +
\mathrm{o}_{\mathbb
{P}}([\Sigma l_i]^{1/2}) \\
& = &\sum_i \sum_j [\varepsilon_{ij}^3 - 3\gamma_{\varepsilon}^2
\varepsilon_{ij}] + \mathrm{o}_{\mathbb{P}}(N^{1/2}) .
\end{eqnarray*}
The conclusion for $\hat\gamma_{\varepsilon}^3$ now follows easily from
the central limit theorem after centering the $\varepsilon_{ij}^3$.
\end{pf*}

\begin{pf*}{Proof of Theorem \ref{thm:6}} For $\hat\gamma_b^4$, we check that
\[
n\hat\gamma_b^4 = \sum_{i=1}^n \frac{1}{l_i(l_i-1)(l_i-2)(l_i-3)}
\sum
_{j\not= k\not= l \not= m} \hat e_{ij} \hat e_{ik} \hat e_{il} \hat
e_{im} .
\]
Again neglecting all terms that contain $z'_{ij}(\beta-\hat\beta)$, we
get, with $v_i = b_i - (\bar{b} + \bar{\varepsilon})$,
\begin{eqnarray*}
n\hat\gamma_b^4 &=& \sum_{i=1}^n \frac{1}{l_i(l_i-1)(l_i-2)(l_i-3)}
\sum
_{j\not= k\not= l \not= m} \varepsilon_{ij} \varepsilon_{ik}
\varepsilon
_{il} \varepsilon_{im} \\
&&{} + 4 \sum_{i=1}^n \frac{v_i}{l_i(l_i-1)(l_i-2)} \sum_{j\not= k\not=
l} \varepsilon_{ij} \varepsilon_{ik} \varepsilon_{il}\\
&&{} + 6 \sum_{i=1}^n \frac{v_i^2}{l_i(l_i-1)} \sum_{j\not= k}
\varepsilon
_{ij} \varepsilon_{ik} + 4\sum_{i=1}^n \frac{v_i^3}{l_i} \sum_j
\varepsilon_{ij} \\
&&{} + \sum_{i=1}^n v_i^4 .
\end{eqnarray*}
Since the first four sums are all $\mathrm{o}_{\mathbb{P}}(n^{1/2})$, we obtain
\[
n^{1/2} [\hat\gamma_b^4 - \gamma_b^4] = n^{-1/2} \sum_{i=1}^n
[v_i^4 -
\gamma_b^4] + \mathrm{o}_{\mathbb{P}}(1) .
\]
The distributional convergence of $\hat\gamma_b^4$ now readily follows
from the central limit theorem after an expansion of the last sum into
\[
n^{-1/2} \sum_{i=1}^n [ b_i^4 - \gamma_b^4 - 4\gamma_b^3 b_i ] +
\mathrm{o}_{\mathbb{P}}(1) .
\]
For $\hat\gamma_{\varepsilon}^4$, we check that
\begin{eqnarray*}
N \hat\gamma_{\varepsilon}^4 &=& -3 \sum_{i=1}^n \frac
{1}{(l_i-1)(l_i-2)(l_i-3)} \sum_{j\not= k\not= l\not= m} \hat e_{ij}
\hat e_{ik} \hat e_{il} \hat e_{im} \\[3pt]
&&{} - 4 \sum_{i=1}^n \sum_{i=1}^n \frac{1}{l_i-1} \sum_{j\not= k}
\hat
e_{ij}^3 \hat e_{ik} + 6 \sum_{i=1}^n \frac{1}{(l_i-1)} \frac
{1}{(l_i-2)} \sum_{j\not= k\not= l} \hat e_{ij}^2 \hat e_{ik} \hat
e_{il} \\[3pt]
&&{} + \sum_{i=1}^n \sum_j \hat e_{ij}^4 .
\end{eqnarray*}
Similarly to the proof of Theorem \ref{thm:5}, it can be shown that
\[
N (\hat\gamma_{\varepsilon}^4 - \gamma_{\varepsilon}^4) = \sum_i
\sum_j
[ \varepsilon_{ij}^4 - \gamma_{\varepsilon}^4 - 4 \gamma
_{\varepsilon}^3 \varepsilon_{ij}] + \mathrm{o}_{\mathbb{P}}(N^{1/2}) ,
\]
from which the conclusion follows.\vspace*{3pt}
\end{pf*}

\begin{pf*}{Proof of Theorem \ref{thm:3}} We first deal with $\hat\gamma
_b^{\ast
3}$. By definition,
\[
\hat\gamma_b^{\ast3} = \frac{1}{n} \sum_{i=1}^n \frac{1}{l_i(l_i
- 1)}
\sum\sum_{j\not= k} \hat e_{ij}^2 \hat e_{ik} .
\]
Our goal will be to use (\ref{eq2.11}) in order to express the last
double sum in terms of $b_i$,~$\varepsilon_{ij}$ and negligible
remainders. Actually, in view of (\ref{eq2.7})--(\ref{eq2.9}), since
the standardizing factor of $\hat\gamma_b^{\ast3}$ is $n^{1/2} =
\mathrm{o}(N^{1/2})$, all terms in the expansion of $\hat\gamma_b^{\ast3}$
containing $z'_{ij}(\beta- \hat\beta)$ are negligible. In other words,
\[
\hat\gamma_b^{\ast3} = \frac{1}{n} \sum_{i=1}^n \frac{1}{l_i(l_i
- 1)}
\sum\sum_{j\not= k} \bigl(b_i + \varepsilon_{ij} - (\bar{b} + \bar
{\varepsilon})\bigr)^2 \bigl(b_i + \varepsilon_{ik} - (\bar{b} + \bar
{\varepsilon
})\bigr) + \mathrm{o}_{\mathbb{P}}(n^{-1/2}) .
\]
After some simple but tedious rearrangements, this becomes
\begin{eqnarray*}
\hat\gamma_b^{\ast3} &=& \frac{1}{n} \sum_{i=1}^n \frac{1}{l_i(l_i-1)}
\biggl[ \sum_{j\not= k} (\varepsilon_{ij}^2 - \gamma_{\varepsilon}^2)
\varepsilon_{ik} + 2b_i \sum_{j\not= k} \varepsilon_{ij}
\varepsilon_{ik} - 2(\bar{b}+\bar{\varepsilon}) \sum_{j\not= k}
\varepsilon_{ij} \varepsilon_{ik}\biggr] \\
&&{} + \frac{1}{n} \sum_{i=1}^n \frac{1}{l_i} \biggl[ b_i \sum_j
(\varepsilon_{ij}^2 - \gamma_{\varepsilon}^2) + 3(b_i^3 - \gamma_b^2)
\sum_j \varepsilon_{ij} \\
&&\hphantom{{} + \frac{1}{n} \sum_{i=1}^n \frac{1}{l_i} \biggl[}{} - (\bar{b} + \bar{\varepsilon}) \sum_j
(\varepsilon
_{ij}^2 - \gamma_{\varepsilon}^2 )- (\bar{b} +\bar{\varepsilon}) 6b_i
\sum_{j} \varepsilon_{ij} \biggr]\\
&&{} - 3(\bar{b}+\bar{\varepsilon}) \frac{1}{n} \sum_{i=1}^n (b_i^2 -
\gamma_b^2) + 2(\bar{b} + \bar{\varepsilon})^3 + \frac{1}{n} \sum
_{i=1}^n b_i^3 - 3\bar{b} \gamma_b^2 + \mathrm{o}_{\mathbb{P}}(n^{-1/2}) .
\end{eqnarray*}
To identify remainders, we note that $n^{1/2}(\bar{b} +
\bar{\varepsilon}) = \mathrm{O}_{\mathbb{P}}(1)$. Also, all summands in the
double and triple sums are centered and independent. Computation of
variances shows that they are all negligible. In summary, we get
\[
n^{1/2} [ \hat\gamma_b^{\ast3} - \gamma_b^3] = n^{-1/2} \sum
_{i=1}^n [
b_i^3 - \gamma_b^3 - 3\gamma_b^2 b_i] + \mathrm{o}_{\mathbb{P}}(1) .
\]
This i.i.d. representation of $\hat\gamma_b^{\ast3}$ is the key tool
for (\ref{eq15}) -- just apply the central limit theorem to the leading
sum.

We will only study $\hat\gamma_{\varepsilon}^{\ast3}$ briefly. First,
by definition,
\[
N \hat\gamma_{\varepsilon}^{\ast3} = \sum_{i=1}^n \sum_{j=1}^{l_i}
\hat e_{ij}^3 - \sum_{i=1}^n \frac{\sum_{j\not= k} \hat e_{ij}^2
\hat
e_{ik}}{l_i - 1} .
\]
To expand the two expressions into leading terms and remainders, recall
that the final standardizing factor in (\ref{eq14}) will be $N^{1/2}$,
which is the same as in (\ref{eq2.9}). We conclude that, under the
conditions of the theorem, terms containing higher orders of
$z'_{ij}(\beta- \hat\beta)$ are negligible. Hence, up to remainders,
\begin{eqnarray*}
\sum_{i=1}^n \sum_{j=1}^{l_i} \hat e_{ij}^3 &=& \sum_{i=1}^n \sum
_{j=1}^{l_i} \bigl(\varepsilon_{ij} - \bar{\varepsilon} + b_i - \bar{b} +
z'_{ij}(\beta- \hat\beta)\bigr)^3\\
&=& \sum_{i=1}^n \sum_{j=1}^{l_i} \varepsilon_{ij}^3 + \sum_{i=1}^n
\sum
_{j=1}^{l_i} 3 \varepsilon_{ij}^2 (-\bar{\varepsilon} + b_i - \bar{b})
+ \sum_{i=1}^n \sum_{j=1}^{l_i} 3\varepsilon_{ij} (-\bar
{\varepsilon} +
b_i - \bar{b})^2 \\
&&{}+ \sum_{i=1}^n \sum_{j=1}^{l_i} (-\bar{\varepsilon} + b_i - \bar
{b})^3 + \sum_{i=1}^n \sum_{j=1}^{l_i} 3 \varepsilon_{ij}^2 z'_{ij}
(\beta- \hat\beta) \\
&&{}+ \sum_{i=1}^{n}\sum_{j=1}^{l_i} 3 (-\bar
{\varepsilon} + b_i - \bar{b})^2 z'_{ij}(\beta- \hat\beta)+ \sum_{i=1}^n \sum_{j=1}^{l_i} 6 \varepsilon_{ij} (- \bar
{\varepsilon} + b_i - \bar{b}) z'_{ij} (\beta- \hat\beta) .
\end{eqnarray*}
A detailed study of these sums yields
\begin{eqnarray*}
&& N^{1/2} (\hat\gamma_{\varepsilon}^{\ast3} - \gamma_{\varepsilon
}^3) \\
&&\quad= N^{-1/2} \sum_{i=1}^n \Biggl\{ 2 \gamma_{\varepsilon}^2 (l_i -
\bar
{l}_n) b_i + \sum_{j=1}^{l_i} \biggl[ \varepsilon_{ij}^3 - \gamma
_{\varepsilon}^3 + 2 b_i (\varepsilon_{ij}^2 - \gamma_{\varepsilon}^2)\\
&&\hphantom{\quad= N^{-1/2} \sum_{i=1}^n \Biggl\{ 2 \gamma_{\varepsilon}^2 (l_i -
\bar
{l}_n) b_i + \sum_{j=1}^{l_i} \biggl[}{}- \gamma_{\varepsilon}^2 \biggl( 1 + \frac{ 2 \bar{l}_n}{l_i} + 2 x'_0
\hat\Sigma^{-1}_n (x_{ij} - \bar{x}_{i\cdot})\biggr)
\varepsilon_{ij}\biggr] \Biggr\} + \mathrm{o}_{\mathbb{P}}(1) .
\end{eqnarray*}
The leading part is a sum of centered independent random variables to
which the central limit theorem may be applied. Its variance satisfies
\begin{eqnarray*}
&& \mathbb{E} [ \varepsilon^3 - \gamma_{\varepsilon}^3]^2 + 4 \gamma_b^2
\mathbb{E} [ \varepsilon^2 - \gamma_{\varepsilon
}^2]^2 + 4 (\gamma_{\varepsilon}^2)^2 \gamma_b^2 d \\
&&\quad{} - 6 \gamma_{\varepsilon}^2 \mathbb{E} \varepsilon^4 + (\gamma
_{\varepsilon}^2)^3 [ 5 + 4c + 4x'_0 \Sigma^{-1} x_0
] + \mathrm{o}(1)\\
&&\quad{} \to\mu_{\varepsilon}^{\ast3} + 4 \gamma_b^2 \bigl( \gamma
_{\varepsilon}^4
- (1-d) (\gamma_{\varepsilon}^2)^2\bigr) ,
\end{eqnarray*}
as desired. This completes the proof of Theorem \ref{thm:3}.
\end{pf*}

\begin{pf*}{Proof of Theorem \ref{thm:4}} The necessary arguments are similar to
those used before and are therefore omitted. Details may be obtained
from the authors.
\end{pf*}
\end{appendix}

\section*{Acknowledgment}
This work was supported by Grant HKBU2030/07p from the Research Grants Council of
Hong Kong.

\printhistory

\end{document}